\documentstyle[11pt,amstex,amssymb]{amsart}
\newtheorem{thm}{Theorem}[section]
\newtheorem{prop}[thm]{Proposition}
\newtheorem{lemma}[thm]{Lemma}
\newtheorem{cor}[thm]{Corollary}
\newtheorem{remark}[thm]{Remark}

\newtheorem{definition}[thm]{Definition}

\newcommand{\proof}{\noindent{\it Proof.}\enspace}

\def\bR{{\Bbb R}}
\def\bN{{\Bbb N}}
\def\bC{{\Bbb C}}
\def\cM{{\cal M}}
\def\cH{{\cal H}}
\def\eps{\varepsilon}
\def\tr{{\rm tr}}
\def\Re{{\rm Re}}
\def\Im{{\rm Im}}
\def\<{\langle}
\def\>{\rangle}
\def\cA{{\cal A}}
\def\1{{\bf 1}}
\def\cT{{\cal T}}
\def\cU{{\cal U}}

\def\min{{\rm min}}

\begin{document}
\title[Free analog of pressure]
{Free analog of pressure and its Legendre transform}
\author[F. Hiai]{Fumio Hiai$\,^1$}
\address{Graduate School of Information Sciences,
Tohoku University, Aoba-ku, Sendai 980-8579, Japan}
\thanks{$^1\,$Supported in part by
Grant-in-Aid for Scientific Research (C)14540198 and by the program
``R\&D support scheme for funding selected IT proposals" of the
Ministry of Public Management, Home Affairs, Posts and
Telecommunications.}

\maketitle

\begin{abstract}
The free analog of the pressure is introduced for multivariate
noncommutative random variables and its Legendre transform is compared with
Voiculescu's microstate free entropy.
\end{abstract}

\section*{Introduction}

The entropy and the pressure are two fundamental ingredients in both
classical and quantum statistical mechanics, in particular, in classical
and quantum lattice systems (see \cite{Is,BR} for example). They are the
dual concepts of each other; more precisely, the entropy function is the
Legendre transform of the pressure function and vice versa under a certain
duality between the state space and the potential space, and an equilibrium
state associated with a potential is usually described by the so-called
variational principle (the equality case of the Legendre transform). The
free entropy introduced by D.~Voiculescu \cite{V1,V2} has played a
central role in free probability theory as the free analog of the
Boltzmann-Gibbs entropy in classical theory. It then would be natural to
consider the free probabilistic analog of the pressure. In \cite{HMP} we
indeed introduced the free pressure of real continuous functions on the
interval $[-R,R]$ and showed its properties like the statistical mechanical
pressure (see Section 1 of this paper).

The aim of the present paper is to introduce the notion of free pressure
for multivariate noncommutative random variables and to investigate it in
connection with the free entropy. (The contents of Sections 2 and 3 were
indeed announced in \cite[\S\S4.4]{Hi}.) In \cite{HMP} we adopted the
Legendre transform of the minus free entropy of probability measures to
define the free pressure of real continuous functions. The idea here is
opposite; we will first introduce the free pressure of noncommutative
multivariables in the so-called microstate approach, and then we will
examine what is the Legendre transform of the free pressure.

In Section 2 we define, given $N\in\bN$ and $R>0$, the free pressure
$\pi_R(h)$ for selfadjoint elements $h$ of the $N$-fold full free product
$C^*$-algebra $\cA_R^{(N)}:=C([-R,R])^{\star N}$ and give its basic
properties. In Sections 3 and 4 we consider the Legendre transform
$\eta_R(\mu)$ of $\pi_R$ for tracial states $\mu$ on $\cA_R^{(N)}$ under
the duality between the selfadjoint elements and the tracial states. For an
$N$-tuple $(a_1,\dots,a_N)$ of selfadjoint noncommutative random variables
in a $W^*$-probability space $(\cM,\tau)$ such that $\|a_i\|\le R$, a
tracial state $\mu_{(a_1,\dots,a_N)}$ on $\cA_R^{(N)}$ can be defined by
$\mu_{(a_1,\dots,a_N)}(h):=\tau(h(a_1,\dots,a_N))$ for $h\in\cA_R^{(N)}$
where $h(a_1,\dots,a_N)$ is the noncommutative ``functional calculus" of
$(a_1,\dots,a_N)$. We then define the free entropy-like quantity
$\eta_R(a_1,\dots,a_N)$ as $\eta_R(\mu_{(a_1,\dots,a_N)})$ and also
$\eta(a_1,\dots,a_N):=\sup_{R>0}\eta_R(a_1,\dots,a_N)$. The properties of
$\eta_R(a_1,\dots,a_N)$ are similar to those of Voiculescu's microstate
free entropy $\chi(a_1,\dots,a_N)$ while they do not generally coincide.
But it is shown that $\eta_R(a_1,\dots,a_N)\ge\chi(a_1,\dots,a_N)$ holds
and equality arises when $a_1,\dots,a_N$ are free. Also, we have
$\eta_R(a_1,a_2)=\chi(a_1,a_2)$ if $a_1+ia_2$ is $R$-diagonal (Section 5).
In Section 6 we slightly modify $\pi_R$ to define the free pressure
$\pi_R^{(2)}(g)$ for selfadjoint elements $g$ of
$\cA_R^{(N)}\otimes_\min\cA_R^{(N)}$ and prove that the quantity
$\tilde\eta(a_1,\dots,a_N)$ induced from $\pi_R^{(2)}$ via Legendre
transform is equal to $\chi(a_1,\dots,a_N)$. In this way, the free entropy
can be understood as the Legendre transform of a certain free probabilistic
pressure. Finally in Section 7 we consider the Gibbs probability measure on
the $N$-fold product of $\bigl\{A\in M_n^{sa}:\|A\|\le R\bigr\}$ associated
with $h_0\in\bigl(\cA_R^{(N)}\bigr)^{sa}$, and we examine the asymptotic
behavior of its Boltzmann-Gibbs entropy as $n\to\infty$ in relation to
$\eta_R(\mu_0)$ of an equilibrium tracial state $\mu_0$ associated with
$h_0$, i.e., a tracial state $\mu_0$ on $\cA_R^{(N)}$ satisfying
$\pi_R(h_0)=-\mu_0(h_0)+\eta_R(\mu_0)$.

\section{Preliminaries}
\setcounter{equation}{0}

Let $(\cM,\tau)$ be a tracial $W^*$-probability space, that is, $\cM$ is a
von Neumann algebra with a faithful normal tracial state $\tau$, and
$\cM^{sa}$ be the set of selfadjoint elements in $\cM$. Let $M_n$ be the
algebra of $n\times n$ complex matrices and $M_n^{sa}$ the set of
selfadjoint matrices in $M_n$. The normalized trace of $A\in M_n$ is
denoted by $\tr_n(A)$ and the operator norm of $A$ by $\|A\|$. In
\cite{V2} D.~Voiculescu introduced the free entropy of an $N$-tuple
$(a_1,\dots,a_N)$ of noncommutative random variables in $\cM^{sa}$ as
follows: For each $R>0$, $\eps>0$ and $n,r\in\bN$ define
\begin{eqnarray*}
&&\Gamma_R(a_1,\dots,a_N;n,r,\eps)
:=\bigl\{(A_1,\dots,A_N)\in(M_n^{sa})^N:\|A_i\|\le R, \\
&&\hskip2cm
|\tr_n(A_{i_1}\cdots A_{i_k})-\tau(a_{i_1}\cdots a_{i_k})|\le\eps,
\ 1\le i_1,\dots,i_k\le N,\ k\le r\bigr\},
\end{eqnarray*}
\begin{eqnarray}\label{F-1.1}
&&\chi_R(a_1,\dots,a_N)
:=\lim_{r\to\infty,\,\eps\to+0}\,\limsup_{n\to\infty}
\biggl({1\over n^2}\log\Lambda_n^{\otimes N}
\bigl(\Gamma_R(a_1,\dots,a_N;n,r,\eps)\bigr) \nonumber\\
&&\hskip10cm+{N\over2}\log n\biggr),
\end{eqnarray}
where $\Lambda_n^{\otimes N}$ denotes the $N$-fold tensor product of the
``Lebesgue" measure $\Lambda_n$ on $M_n^{sa}$:
$$
d\Lambda_n(A):=2^{n(n-1)/2}\prod_{i=1}^ndA_{ii}
\prod_{i<j}d(\Re\,A_{ij})\,d(\Im\,A_{ij})
$$
(the constant $2^{n(n-1)/2}$ comes from the isometric isomorphism between
$M_n^{sa}$ and $\bR^{n^2}$). Then the {\it free entropy} of
$(a_1,\dots,a_N)$ is
$$
\chi(a_1,\dots,a_N):=\sup_{R>0}\chi_R(a_1,\dots,a_N).
$$
The definition being based on the microstate (or matricial) approximation,
this\break $\chi(a_1,\dots,a_N)$ is sometimes called the {\it microstate
free entropy} in contrast to another Voiculescu's free entropy
$\chi^*(a_1,\dots,a_N)$ in the microstate-free approach in \cite{V5}.

In the case of a single variable $a\in\cM^{sa}$ whose distribution measure
(with respect to $\tau$) is $\mu$, $\chi(a)$ coincides with the free
entropy $\Sigma(\mu):=\iint\log|x-y|\,d\mu(x)\,d\mu(y)$ of $\mu$ introduced
in \cite{V1} up to an additive constant:
$$
\chi(a)=\Sigma(\mu)+{1\over2}\log2\pi+{3\over4}.
$$

With $R>0$ fixed, let $C_\bR([-R,R])$ be the Banach space of real continuous
functions on $[-R,R]$ with sup-norm $\|\cdot\|$, and $\cM([-R,R])$ be the
set of probability measures on $[-R,R]$. Consider the dual pairing
$$
\mu(h):=\int h\,d\mu
\quad\mbox{for $h\in C_\bR([-R,R])$, $\mu\in\cM([-R,R])$}.
$$
The free entropy $\chi(\mu):=\Sigma(\mu)+{1\over2}\log2\pi+{3\over4}$ being
strictly concave and weakly upper semicontinuous on $\cM([-R,R])$ (see
\cite[5.3.2]{HP1} for example), it is natural to consider the {\it Legendre
transform} of $-\chi(\mu)$ as follows:
\begin{equation}\label{F-1.2}
\pi_R(h):=\sup\bigl\{-\mu(h)+\chi(\mu):\mu\in\cM([-R,R])\bigr\}
\quad\mbox{for $h\in C_\bR([-R,R])$}.
\end{equation}
This $\pi_R(h)$ is the {\it free pressure} of $h$ discussed in \cite{HMP}.
A fundamental result in the theory of weighted potentials (\cite[I.1.3 and
I.3.1]{ST}) tells us that for each $h\in C_\bR([-R,R])$ there exists a
unique $\sigma^h\in\cM([-R,R])$, called the {\it equilibrium measure}
associated with $h$, such that
$$
-\sigma^h(h)+\chi(\sigma^h)=\pi_R(h).
$$
The last equality characterizing $\sigma^h$ is a kind of the variational
principle, a fundamental notion in statistical mechanics. Let us assemble
some properties of the free pressure obtained in \cite{HMP} in the following:

\begin{itemize}
\item[(I)] The function $\pi_R(h)$ on $C_\bR([-R,R])$ is convex and
decreasing, i.e., $\pi_R(h_1)\ge\pi_R(h_2)$ if $h_1\le h_2$. Moreover,
$$
|\pi_R(h_1)-\pi_R(h_2)|\le\|h_1-h_2\|
\quad\mbox{for all $h_1,h_2\in C_\bR([-R,R])$}.
$$
\item[(II)] Conversely, $\chi(\mu)$ is the (minus) Legendre transform of
$\pi_R(h)$, that is,
\begin{equation}\label{F-1.3}
\chi(\mu)=\inf\bigl\{\mu(h)+\pi_R(h):h\in C_\bR([-R,R])\bigr\}
\quad\mbox{for every $\mu\in\cM([-R,R])$},
\end{equation}
\item[(III)] For every $h\in C_\bR([-R,R])$, $\pi_R(h)$ is expressed as
\begin{equation}\label{F-1.4}
\pi_R(h)=\lim_{n\to\infty}\Biggl({1\over n^2}\log
\int_{(M_n^{sa})_R}\exp\bigl(-n^2\tr_n(h(A))\bigr)\,d\Lambda_n(A)
+{1\over2}\log n\Biggr),
\end{equation}
where $(M_n^{sa})_R:=\bigl\{A\in M_n^{sa}:\|A\|\le R\bigr\}$ and $h(A)$ is
defined via functional calculus.
\item[(IV)] For each $h\in C_\bR([-R,R])$, define a probability measure
$\lambda_{R,n}^h$ on $(M_n^{sa})_R$ by
$$
d\lambda_{R,n}^h(A):={1\over Z_{R,n}^h}\exp\bigl(-n^2\tr_n(h(A))\bigr)
\,\chi_{\{\|A\|\le R\}}(A)\,d\Lambda_n(A)
$$
with the normalization constant
$$
Z_{R,n}^h:=\int_{(M_n^{sa})_R}\exp\bigl(-n^2
\tr_n(h(A))\bigr)\,d\Lambda_n(A),
$$
and let $S(\lambda_{R,n}^h)$ be the Boltzmann-Gibbs entropy of
$\lambda_{R,n}^h$. Then
$$
\chi(\sigma^h)=\lim_{n\to\infty}\biggl(
{1\over n^2}S(\lambda_{R,n}^h)+{1\over2}\log n\biggr).
$$
\end{itemize}

\section{Free pressure $\pi_R$ for multivariables}
\setcounter{equation}{0}

For each $N\in\bN$ and $R>0$ we write $\cA_R^{(N)}$ for the $N$-fold full
(universal) free product $C^*$-algebra of $C([-R,R])$, i.e.,
$\cA_R^{(N)}:=C([-R,R])^{\star N}$. One can describe $\cA_R^{(N)}$ in a bit
more constructive way as below. Consider the (algebraic) noncommutative
polynomial $*$-algebra $\bC\<X_1,\dots,X_N\>$ with noncommuting
indeterminates $X_1,\dots,X_N$, where the $*$-operation is given by
$X_i^*=X_i$. Let $\bC\<X_1,\dots,X_N\>^{sa}$ be the
space of selfadjoint polynomials $p=p^*$ in $\bC\<X_1,\dots,X_N\>$.
Obviously, $p(A_1,\dots,A_N)\in M_n^{sa}$ whenever
$p\in\bC\<X_1,\dots,X_N\>^{sa}$ and $A_1,\dots,A_N\in M_n^{sa}$. One can
immediately see that a $C^*$-norm $\|\cdot\|_R$ is defined on
$\bC\<X_1,\dots,X_N\>$ by
$$
\|p\|_R:=\sup\bigl\{\|p(A_1,\dots,A_N)\|:A_1,\dots,A_N\in M_n^{sa},
\ \|A_i\|\le R,\ n\in\bN\bigr\}.
$$

\begin{prop}\label{P-2.1}
The $C^*$-completion of $\bC\<X_1,\dots,X_N\>$ with respect to $\|\cdot\|_R$
is isomorphic to $\cA_R^{(N)}$ by the isomorphism mapping $X_i$ to
$f_0(t)=t$ in the $i$th copy of $C([-R,R])$ of $\cA_R^{(N)}$ for
$1\le i\le N$.
\end{prop}

\proof
Let ${\frak A}_R^{(N)}$ be the $C^*$-completion of $\bC\<X_1,\dots,X_N\>$
with respect to $\|\cdot\|_R$. Since any $*$-homomorphism of $C([-R,R])$
into $B(\cH)$ is determined by a selfadjoint $a\in B(\cH)$ with
$\|a\|\le R$ (the image of $f_0$), what we have to prove is that, for each
selfadjoint $a_1,\dots,a_N\in B(\cH)$ with $\|a_i\|\le R$, there exists
a $*$-homomorphism $\Phi:{\frak A}_R^{(N)}\to B(\cH)$ such that
$\Phi(X_i)=a_i$, $1\le i\le N$. To do so, it suffices to show that
$$
\|p(a_1,\dots,a_N)\|\le\|p\|_R
\quad\mbox{for all $p\in\bC\<X_1,\dots,X_N\>$}.
$$
But this is easy to show. In fact, choose a net of finite-dimensional
projections $e_k$ in $B(\cH)$ with $e_k\nearrow\1$. Then, the above
inequality follows because $\|p(e_ka_1e_k,\dots,e_ka_Ne_k)\|\le\|p\|_R$ and
$p(e_ka_1e_k,\dots,e_ka_Ne_k)\to p(a_1,\dots,a_N)$ strongly.\qed

\bigskip
In view of the above proposition, we consider $\cA_R^{(N)}$ as the
$C^*$-completion of $\bC\<X_1,\dots,X_N\>$ with respect to $\|\cdot\|_R$ and
write $\|\cdot\|_R$ for the norm of $\cA_R^{(N)}$ as well. For each
selfadjoint $a_1,\dots,a_N\in B(\cH)$ with $\|a_i\|\le R$, we can
extend $p\in\bC\<X_1,\dots,X_N\>\mapsto p(a_1,\dots,a_N)\in B(\cH)$ to
a $*$-homomorphism $\cA_R^{(N)}\to B(\cH)$, written as
$$
h\in\cA_R^{(N)}\mapsto h(a_1,\dots,a_N)\in B(\cH),
$$
which is regarded as the ``continuous functional calculus" of the
noncommutative $N$-tuple $(a_1,\dots,a_N)$. In particular, for a single
selfadjoint $a$ with $\|a\|\le R$, this reduces to the usual functional
calculus $h(a)$ for $h\in C([-R,R])$ ($=\cA_R^{(1)}$). There are natural
operations on the class of $C^*$-algebras $\cA_R^{(N)}$ ($N\in\bN$, $R>0$).
\begin{itemize}
\item[(a)] When $R_1<R$, the restriction $f\mapsto f|_{[-R_1,R_1]}$,
$f\in C([-R,R])$, induces a $*$-homomorphism
$r_{R,R_1}:\cA_R^{(N)}\to\cA_{R_1}^{(N)}$. Note that $r_{R,R_1}(p)=p$ for
$p\in\bC\<X_1,\dots,X_N\>$.
\item[(b)] For each $R,R_1>0$ the dilation $f\mapsto f((R/R_1)\,\cdot)$
induces a $*$-isomorphism $\rho_{R,R_1}:\cA_R^{(N)}\to\cA_{R_1}^{(N)}$. For
a polynomial $p\in\bC\<X_1,\dots,X_N\>$ this is written as
$$
\qquad\rho_{R,R_1}(p)(X_1,\dots,X_N)=p((R/R_1)X_1,\dots,(R/R_1)X_N).
$$
In particular, the $C^*$-algebra $\cA_R^{(N)}$ is independent (up to a
$*$-isomorphism) of the choice of $R>0$.
\item[(c)] When $1\le L<N$, $\cA_R^{(L)}$ and $\cA_R^{(N-L)}$ are considered
as $C^*$-subalgebras of $\cA_R^{(N)}$ under the identification
$\cA_R^{(L)}\star\cA_R^{(N-L)}\cong\cA_R^{(N)}$. For $h_1\in\cA_R^{(L)}$ and
$h_2\in\cA_R^{(N-L)}$ we simply write $h_1+h_2$ for the sum of $h_1,h_2$ in
$\cA_R^{(N)}$ under this identification. For polynomials
$p_1\in\bC\<X_1,\dots,X_L\>$ and $p_2\in\bC\<X_{L+1},\dots,X_N\>$ this is
the usual polynomial sum $p_1+p_2\in\bC\<X_1,\dots,X_N\>$.
\end{itemize}

Let $\bigl(\cA_R^{(N)}\bigr)^{sa}$ be the selfadjoint part of $\cA_R^{(N)}$.
Note that $h(A_1,\dots,A_N)\in M_n^{sa}$ whenever
$h\in\bigl(\cA_R^{(N)}\bigr)^{sa}$ and $A_1,\dots,A_N\in M_n^{sa}$ with
$\|A_i\|\le R$. It is immediate to see that the function $h(A_1,\dots,A_N)$
on $(M_n^{sa})_R^N$ is continuous, where $(M_n^{sa})_R^N$ is the $N$-fold
product of $(M_n^{sa})_R$. The following definition is the direct
multivariate extension of the formula \eqref{F-1.4} (see also \eqref{F-1.1}):

\begin{definition}\label{D-2.2}{\rm
For each $R>0$ and $h\in\bigl(\cA_R^{(N)}\bigr)^{sa}$ define
\begin{eqnarray}\label{F-2.1}
&&\pi_R(h):=\limsup_{n\to\infty}
\Biggl({1\over n^2}\log\int_{(M_n^{sa})_R^N}
\exp\bigl(-n^2\tr_n(h(A_1,\dots,A_N))\bigr)
\,d\Lambda_n^{\otimes N}(A_1,\dots,A_N) \nonumber\\
&&\hskip10cm+{N\over2}\log n\Biggr).
\end{eqnarray}
We call this $\pi_R(h)$ the {\it free pressure} of $h$ (with parameter
$R>0$).
}\end{definition}

Basic properties of free entropy $\pi_R(h)$ are summarized in the following.
This and the above definition extend the properties (I) and (III) in Section
1 to noncommutative multivariables.

\begin{prop}\label{P-2.3} \
\begin{itemize}
\item[(1)] $\pi_R(h)$ is convex on $\bigl(\cA_R^{(N)}\bigr)^{sa}$.
\item[(2)] If $h_1,h_2\in\bigl(\cA_R^{(N)}\bigr)^{sa}$ and $h_1\le h_2$,
then $\pi_R(h_1)\ge\pi_R(h_2)$.
\item[(3)] $|\pi_R(h_1)-\pi_R(h_2)|\le\|h_1-h_2\|_R$ for all
$h_1,h_2\in\bigl(\cA_R^{(N)}\bigr)^{sa}$.
\item[(4)] If $R>R_1>0$ and $h\in\bigl(\cA_R^{(N)}\bigr)^{sa}$, then
$\pi_{R_1}(r_{R,R_1}(h))\le\pi_R(h)$.
\item[(5)] If $R,R_1>0$ and $h\in\bigl(\cA_R^{(N)}\bigr)^{sa}$, then
$\pi_{R_1}(\rho_{R,R_1}(h))=\pi_R(h)+N\log(R_1/R)$.
\item[(6)] Let $h_1\in\bigl(\cA_R^{(L)}\bigr)^{sa}$ and
$h_2\in\bigl(\cA_R^{(N-L)}\bigr)^{sa}$ where $1\le L<N$, and $h_1+h_2$ be
given as in (c) above. Then $\pi_R(h_1+h_2)\le\pi_R(h_1)+\pi_R(h_2)$. In
particular when $L=1$, $\pi_R(h_1+h_2)=\pi_R(h_1)+\pi_R(h_2)$.
\end{itemize}
\end{prop}

\proof
(1)\enspace Let $h_1,h_2\in\bigl(\cA_R^{(N)}\bigr)^{sa}$ and $0<\alpha<1$.
For
$A_1,\dots,A_N\in(M_n^{sa})_R$ set
\begin{eqnarray*}
f_k(A_1,\dots,A_N)&:=&\exp\bigl(-n^2\tr_n(h_k(A_1,\dots,A_N))\bigr),
\qquad k=1,2, \\
f(A_1,\dots,A_N)&:=&\exp\bigl(-n^2\tr_n(
(\alpha h_1+(1-\alpha)h_2)(A_1,\dots,A_N))\bigr).
\end{eqnarray*}
Since the H\"older inequality gives
\begin{eqnarray*}
\int_{(M_n^{sa})_R^N}f\,d\Lambda_n^{\otimes N}
&=&\int_{(M_n^{sa})_R^N}f_1^\alpha f_2^{1-\alpha}
\,d\Lambda_n^{\otimes N} \\
&\le&\Biggl(\int_{(M_n^{sa})_R^N}f_1\,d\Lambda_n^{\otimes N}\Biggr)^\alpha
\Biggl(\int_{(M_n^{sa})_R^N}f_2\,d\Lambda_n^{\otimes N}\Biggr)^{1-\alpha},
\end{eqnarray*}
we have
\begin{eqnarray*}
&&{1\over n^2}\log\int_{(M_n^{sa})_R^N}f\,d\Lambda_n^{\otimes N}
+{N\over2}\log n \\
&&\qquad\le\alpha\Biggl({1\over n^2}\log
\int_{(M_n^{sa})_R^N}f_1\,d\Lambda_n^{\otimes N}+{N\over2}\log n\Biggr) \\
&&\qquad\qquad+(1-\alpha)\Biggl({1\over n^2}\log
\int_{(M_n^{sa})_R^N}f_2\,d\Lambda_n^{\otimes N}+{N\over2}\log n\Biggr),
\end{eqnarray*}
which implies that
$$
\pi_R(\alpha h_1+(1-\alpha)h_2)\le\alpha\pi_R(h_1)+(1-\alpha)\pi_R(h_2).
$$

(2) is obvious because the assumption implies that
$h_1(A_1,\dots,A_N)\le h_2(A_1,\dots,A_N)$ for all $A_i\in(M_n^{sa})_R$.

(3)\enspace By the definition of $\|\cdot\|_R$, for any $A_1,\dots,A_N\in
(M_n^{sa})_R$ we get
$$
|\tr_n(h_1(A_1,\dots,A_N))-\tr_n(h_2(A_1,\dots,A_N))|
\le\|h_1-h_2\|_R
$$
so that
\begin{eqnarray*}
&&\exp\bigl(-n^2\tr_n(h_1(A_1,\dots,A_N))\bigr)
\exp\bigl(-n^2\|h_1-h_2\|_R\bigr) \\
&&\qquad\le\exp\bigl(-n^2\tr_n(h_2(A_1,\dots,A_N))\bigr) \\
&&\qquad\le\exp\bigl(-n^2\tr_n(h_1(A_1,\dots,A_N))\bigr)
\exp\bigl(n^2\|h_1-h_2\|_R\bigr).
\end{eqnarray*}
This immediately implies that
$$
\pi_R(h_1)-\|h_1-h_2\|_R\le\pi_R(h_2)\le\pi_R(h_1)+\|h_1-h_2\|_R.
$$

(4) is obvious because $r_{R,R_1}(h)(A_1,\dots,A_N)=h(A_1,\dots,A_N)$ for
$A_i\in(M_n^{sa})_{R_1}$.

(5)\enspace For $h_1:=\rho_{R,R_1}(h)$ and $\alpha:=R_1/R$, since
$h_1(A_1,\dots,A_N)=h(\alpha^{-1}A_1,\dots,\alpha^{-1}A_N)$ for
$A_i\in(M_n^{sa})_{R_1}$, we get
\begin{eqnarray*}
&&\int_{(M_n^{sa})_{R_1}^N}\exp\bigl(-n^2
\tr_n(h_1(A_1,\dots,A_N))\bigr)\,d\Lambda_n^{\otimes N} \\
&&\qquad=\int_{(M_n^{sa})_{R_1}^N}\exp\bigl(-n^2\tr_n
(h(\alpha^{-1}A_1,\dots,\alpha^{-1}A_N))\bigr)\,d\Lambda_n^{\otimes N} \\
&&\qquad=\alpha^{n^2N}\int_{(M_n^{sa})_R^N}\exp\bigl(-n^2
\tr_n(h(A_1,\dots,A_N))\bigr)\,d\Lambda_n^{\otimes N}
\end{eqnarray*}
thanks to the trivial formula $d\Lambda_n(\alpha
A)/d\Lambda_n(A)=\alpha^{n^2}$. Hence $\pi_{R_1}(h_1)=\pi_R(h)+N\log\alpha$.

(6)\enspace Since
$(h_1+h_2)(A_1,\dots,A_N)=h_1(A_1,\dots,A_L)+h_2(A_{L+1},\dots,A_N)$ for
$A_i\in(M_n^{sa})_R$, we get
\begin{eqnarray*}
&&\int_{(M_n^{sa})_R^N}\exp\bigl(-n^2
\tr_n((h_1+h_2)(A_1,\dots,A_N))\bigr)\,d\Lambda_n^{\otimes N} \\
&&\qquad=\int_{(M_n^{sa})_R^L}\exp\bigl(-n^2
\tr_n(h_1(A_1,\dots,A_L))\bigr)\,d\Lambda_n^{\otimes L} \\
&&\qquad\qquad\times\int_{(M_n^{sa})_R^{N-L}}\exp\bigl(-n^2
\tr_n(h_2(A_{L+1},\dots,A_N))\bigr)\,d\Lambda_n^{\otimes{N-L}}
\end{eqnarray*}
and hence
\begin{eqnarray*}
&&{1\over n^2}\log\int_{(M_n^{sa})_R^N}\exp\bigl(-n^2
\tr_n((h_1+h_2)(A_1,\dots,A_N))\bigr)\,d\Lambda_n^{\otimes N}
+{N\over2}\log n \\
&&\quad={1\over n^2}\log\int_{(M_n^{sa})_R^L}\exp\bigl(-n^2
\tr_n(h_1(A_1,\dots,A_L))\bigr)\,d\Lambda_n^{\otimes L}
+{L\over2}\log n \\
&&\qquad+{1\over n^2}\log\int_{(M_n^{sa})_R^{N-L}}\exp\bigl(-n^2
\tr_n(h_2(A_{L+1},\dots,A_N))\bigr)\,d\Lambda_n^{\otimes{N-L}}
+{N-L\over2}\log n,
\end{eqnarray*}
which gives the required inequality. Also, the above equality together with
the formula \eqref{F-1.4} (with limit) gives the last assertion.
\qed

\section{Legendre transform $\eta_R$ of free pressure}
\setcounter{equation}{0}

We first define the (minus) Legendre transform of the free pressure $\pi_R$
in the purely algebraic situation.

\begin{definition}\label{D-3.1}{\rm
Let $\Sigma\<X_1,\dots,X_N\>$ denote the set of selfadjoint linear
functionals $\mu$ on $\bC\<X_1,\dots,X_N\>$ such that $\mu(\1)=1$, where
the selfadjointness of $\mu$ means $\mu(p^*)=\overline{\mu(p)}$ for
$p\in\bC\<X_1,\dots,X_N\>$. For each $R>0$ and
$\mu\in\Sigma\<X_1,\dots,X_N\>$ define
$$
\eta_R(\mu):=\inf\bigl\{\mu(p)+\pi_R(p):
p\in\bC\<X_1,\dots,X_N\>^{sa}\bigr\}.
$$
}\end{definition}

The following is obvious by definition.

\begin{prop}\label{P-3.2}
For each $R>0$, $\eta_R(\mu)$ is a concave function on
$\Sigma\<X_1,\dots,X_N\>$, and it is upper semicontinuous in the sense that
if $\mu,\mu_k\in\Sigma\<X_1,\dots,X_N\>$ for $k\in\bN$ and
$\mu_k(p)\to\mu(p)$ as $k\to\infty$ for all $p\in\bC\<X_1,\dots,X_N\>$, then
$$
\eta_R(\mu)\ge\limsup_{k\to\infty}\eta_R(\mu_k).
$$
\end{prop}

Let $\cT\bigl(\cA_R^{(N)}\bigr)$ denote the set of all tracial states on
$\cA_R^{(N)}$. As is easily seen, $\mu\in\Sigma\<X_1,\dots,X_N\>$ can
(uniquely) extend to an element of $\cT\bigl(\cA_R^{(N)}\bigr)$ (in this
case we also write $\mu\in\cT\bigl(\cA_R^{(N)}\bigr)$) if and only if the
following three conditions are satisfied:
\begin{itemize}
\item[(a)] $\mu\ge0$ in the sense that $\mu(p^*p)\ge0$ for all
$p\in\bC\<X_1,\dots,X_n\>$,
\item[(b)] $|\mu(p)|\le\|p\|_R$ for all $p\in\bC\<X_1,\dots,X_N\>$,
\item[(c)] $\mu$ is tracial in the sense that $\mu(p_1p_2)=\mu(p_2p_1)$ for
all $p_1,p_2\in\bC\<X_1,\dots,X_N\>$.
\end{itemize}

\begin{lemma}\label{L-3.3}
If $\mu\in\Sigma\<X_1,\dots,X_N\>$ and $\mu\notin\cT\bigl(\cA_R^{(N)}\bigr)$
(i.e., $\mu$ does not extend to an element of $\cT\bigl(\cA_R^{(N)}\bigr)$),
then $\eta_R(\mu)=-\infty$.
\end{lemma}

\proof
Assume that $\mu\in\Sigma\<X_1,\dots,X_N\>$ and $\eta_R(\mu)>-\infty$,
and we prove that the above conditions (a)--(c) hold.

(a)\enspace Suppose $\mu(p^*p)<0$ for some $p\in\bC\<X_1,\dots,X_N\>$; then
for each $\alpha>0$ we get $\pi_R(\alpha p^*p)\le\pi_R(0)$ by Proposition
\ref{P-2.3}\,(2). Therefore,
$$
\eta_R(\mu)\le\mu(\alpha p^*p)+\pi_R(\alpha p^*p)
\longrightarrow-\infty\quad\mbox{as $\alpha\to\infty$},
$$
a contradiction.

(b)\enspace Suppose $|\mu(p)|>\|p\|_R$ for some
$p\in\bC\<X_1,\dots,X_N\>^{sa}$. We may suppose $\mu(p)<-\|p\|_R$. Since
Proposition \ref{P-2.3}\,(3) gives $\pi_R(\alpha p)\le\|\alpha
p\|_R+\pi_R(0)$ for $\alpha>0$, we get
\begin{eqnarray*}
\eta_R(\mu)&\le&\mu(\alpha p)+\pi_R(\alpha p) \\
&\le&\alpha\bigl(\mu(p)+\|p\|_R\bigr)+\pi_R(0)\longrightarrow-\infty
\quad\mbox{as $\alpha\to\infty$},
\end{eqnarray*}
a contradiction. Hence we have shown that $|\mu(p)|\le\|p\|_R$ for all
$p\in\bC\<X_1,\dots,X_N\>^{sa}$. Thanks to this and (a) already proven, it
follows that $\mu$ extends to a bounded positive linear functional on the
$C^*$-algebra $\cA_R^{(N)}$. Since $\mu(\1)=1$, the extended $\mu$ has norm
one.

(c)\enspace It suffices to show that $\mu(i(p_1p_2-p_2p_1))=0$ for all
$p_1,p_2\in\bC\<X_1,\dots,X_N\>^{sa}$. Suppose $\mu(i(p_1p_2-p_2p_1))<0$
for some $p_1,p_2\in\bC\<X_1,\dots,X_n\>^{sa}$. For $\alpha>0$ we get
$\pi_R(i\alpha(p_1p_2-p_2p_1))=\pi_R(0)$ because
$\tr_n(i\alpha(p_1p_2-p_2p_1)(A_1,\dots,A_N))=0$ for any $A_i\in M_n^{sa}$.
Therefore,
\begin{eqnarray*}
\eta_R(\mu)&\le&\mu(i\alpha(p_1p_2-p_2p_1))
+\pi_R(i\alpha(p_1p_2-p_2p_1)) \\
&=&\alpha\mu(i(p_1p_2-p_2p_1))+\pi_R(0)\longrightarrow-\infty
\quad\mbox{as $\alpha\to\infty$},
\end{eqnarray*}
a contradiction.\qed

\bigskip
The above lemma says that the essential domain
$\bigl\{\mu\in\Sigma\<X_1,\dots,X_N\>:\eta_R(\mu)>-\infty\bigr\}$ is
included in $\cT\bigl(\cA_R^{(N)}\bigr)$. In the next theorem we show that
$\pi_R(h)$ for $h\in\bigl(\cA_R^{(N)}\bigr)^{sa}$ and $\eta_R(\mu)$ for
$\mu\in\cT\bigl(\cA_R^{(N)}\bigr)$ are the Legendre transforms of each other
with respect to the Banach space duality $\mu(h)$ for
$h\in\bigl(\cA_R^{(N)}\bigr)^{sa}$ and
$\mu\in\bigl(\cA_R^{(N)}\bigr)^{*,sa}$, the selfadjoint part of
$\bigl(\cA_R^{(N)}\bigr)^*$. In this way, the Legendre transform formulas
\eqref{F-1.2} and \eqref{F-1.3} are extended to the noncommutative
multivariate setting, where $C_\bR([-R,R])$ and $\cM([-R,R])$ are replaced
by $\bigl(\cA_R^{(N)}\bigr)^{sa}$ and $\cT\bigl(\cA_R^{(N)}\bigr)$,
respectively.

\begin{thm}\label{T-3.4} \
\begin{itemize}
\item[(1)] For every $\mu\in\cT\bigl(\cA_R^{(N)}\bigr)$,
$$
\eta_R(\mu)=\inf\bigl\{\mu(h)+\pi_R(h):
h\in\bigl(\cA_R^{(N)}\bigr)^{sa}\bigr\}.
$$
Hence, $\eta_R(\mu)$ is concave and weakly* upper semicontinuous on
$\cT\bigl(\cA_R^{(N)}\bigr)$.
\item[(2)] For every $h\in\bigl(\cA_R^{(N)}\bigr)^{sa}$,
$$
\pi_R(h)=\sup\bigl\{-\mu(h)+\eta_R(\mu):
\mu\in\cT\bigl(\cA_R^{(N)}\bigr)\bigr\}.
$$
\end{itemize}
\end{thm}

\proof
(1) is obvious from Definition \ref{D-3.1} because
$\bC\<X_1,\dots,X_N\>^{sa}$ is dense in $\bigl(\cA_R^{(N)}\bigr)^{sa}$.

(2)\enspace The Legendre transform of $\pi_R(h)$ on
$\bigl(\cA_R^{(N)}\bigr)^{sa}$ with respect to the Banach space duality
between $\bigl(\cA_R^{(N)}\bigr)^{sa}$ and $\bigl(\cA_R^{(N)}\bigr)^{*,sa}$
is
\begin{eqnarray*}
(\pi_R)^*(\mu)&:=&\sup\bigl\{-\mu(h)-\pi_R(h):
h\in\bigl(\cA_R^{(N)}\bigr)^{sa}\bigr\} \\
&=&-\inf\bigl\{\mu(h)+\pi_R(h):h\in\bigl(\cA_R^{(N)}\bigr)^{sa}\bigr\}
\quad\mbox{for $\mu\in\bigl(\cA_R^{(N)}\bigr)^{*,sa}$}.
\end{eqnarray*}
Suppose $\mu\in\bigl(\cA_R^{(N)}\bigr)^{*,sa}$ and $\mu(\1)\ne1$; then for
$\alpha\in\bR$ we get
$$
\mu(\alpha\1)+\pi_R(\alpha\1)=\alpha(\mu(\1)-1)+\pi_R(0)
$$
thanks to $\pi_R(\alpha\1)=-\alpha+\pi_R(0)$. Hence
$(\pi_R)^*(\mu)=+\infty$. This and Lemma \ref{L-3.3} imply that
$$
(\pi_R)^*(\mu)=\cases
-\eta_R(\mu) & \text{if $\mu\in\cT\bigl(\cA_R^{(N)}\bigr)$}, \\
+\infty & \text{if $\mu\in\bigl(\cA_R^{(N)}\bigr)^{*,sa}
\setminus\cT\bigl(\cA_R^{(N)}\bigr)$}.
\endcases
$$
Since $\pi_R(h)$ is convex and continuous on $\bigl(\cA_R^{(N)}\bigr)^{sa}$,
it is the Legendre transform of $(\pi_R)^*(\mu)$ on
$\bigl(\cA_R^{(N)}\bigr)^{*,sa}$ so that
\begin{eqnarray*}
\pi_R(h)
&=&\sup\bigl\{-\mu(h)-(\pi_R)^*(\mu):
\mu\in\bigl(\cA_R^{(N)}\bigr)^{*,sa}\bigr\} \\
&=&\sup\bigl\{-\mu(h)+\eta_R(\mu):\mu\in\cT\bigl(\cA_R^{(N)}\bigr)\bigr\},
\end{eqnarray*}
as desired.\qed

\bigskip
Since $\cT\bigl(\cA_R^{(N)}\bigr)$ is weakly* compact, for each
$h_0\in\bigl(\cA_R^{(N)}\bigr)^{sa}$ there exists a
$\mu_0\in\cT\bigl(\cA_R^{(N)}\bigr)$ such that
$$
\pi_R(h_0)=-\mu_0(h_0)+\eta_R(\mu_0).
$$
This equality condition is a kind of variational principle, so we call such
$\mu_0$ an {\it equilibrium tracial state} associated with $h_0$. We see by
the above proof that such a $\mu_0\in\cT\bigl(\cA_R^{(N)}\bigr)$ is
characterized as an element of $\bigl(\cA_R^{(N)}\bigr)^{*,sa}$ which
supports at $h_0$ the function $\pi_R(h)$ on $\bigl(\cA_R^{(N)}\bigr)^{sa}$.
Thus, the general theory of conjugate (or Legendre) functions says (see
\cite[I.5.3]{ET} for example) that the uniqueness of an equilibrium tracial
state associated with
$h_0\in\cA_R^{sa}$ is equivalent to the differentiability of $\pi_R(h)$ at
$h_0$; or equivalently,
$$
\lim_{t\to0}{\pi_R(h_0+tp)-\pi_R(h_0)\over t}
$$
exists for every $p\in\bC\<X_1,\dots,X_N\>^{sa}$ (the above limit then
equals $\mu_0(p)$). In the single variable case ($N=1$), as mentioned in
Section 1, there is a unique equilibrium measure $\sigma^{h_0}$ associated
with each $h_0\in C_\bR([-R,R])$ ($=\bigl(\cA_R^{(1)}\bigr)^{sa}$), so that
we have
$$
\lim_{t\to0}{\pi_R(h_0+th)-\pi_R(h_0)\over t}=\sigma^{h_0}(h)
\quad\mbox{for every $h\in C_\bR([-R,R])$}.
$$
However, in the multivariate case, it seems quite a difficult problem to
show the differentiability of $\pi_R(h)$ (even at $h_0=0$ for instance). It
might be worth noting a general result that the function $\pi_R$ is
differentiable at points in a dense $G_\delta$ subset of
$\bigl(\cA_R^{(N)}\bigr)^{sa}$. Indeed, this is true for any Lipschitz
continuous convex function on a separable Banach space (see the proof of
\cite[V.9.8]{DS}).

\section{Quantity $\eta(a_1,\dots,a_N)$ for multivariables}
\setcounter{equation}{0}

We apply the quantity $\eta_R(\mu)$ for
$\mu\in\Sigma\<X_1,\dots,X_N\>$ to introduce a free entropy-like quantity
for noncommutative multivariables in a tracial $W^*$-probability space
$(\cM,\tau)$.

\begin{definition}\label{D-4.1}{\rm
For each $a_1,\dots,a_N\in\cM^{sa}$ define
$\mu_{(a_1,\dots,a_N)}\in\Sigma\<X_1,\dots,X_N\>$ by
$$
\mu_{(a_1,\dots,a_N)}(p):=\tau(p(a_1,\dots,a_N)).
$$
We further define
$$
\eta_R(a_1,\dots,a_N):=\eta_R(\mu_{(a_1,\dots,a_N)});
$$
namely,
$$
\eta_R(a_1,\dots,a_N)
=\inf\bigl\{\tau(p(a_1,\dots,a_N))+\pi_R(p):
p\in\bC\<X_1,\dots,X_N\>^{sa}\bigr\},
$$
and
$$
\eta(a_1,\dots,a_N):=\sup_{R>0}\eta_R(a_1,\dots,a_N).
$$
}\end{definition}

If $(a_1,\dots,a_N)\in\cM^{sa}$ and $R\ge\max_i\|a_i\|$, then we have
$\mu:=\mu_{(a_1,\dots,a_N)}\in\cT\bigl(\cA_R^{(N)}\bigr)$ because of the the
functional calculus $h\in\cA_R^{(N)}\mapsto h(a_1,\dots,a_N)\in\cM$. In this
case, $\eta_R(a_1,\dots,a_N)$ is also given as
$$
\eta_R(a_1,\dots,a_N)
=\inf\bigl\{\tau(h(a_1,\dots,a_N))+\pi_R(h):
h\in\bigl(\cA_R^{(N)}\bigr)^{sa}\bigr\},
$$
and we further see that $\{a_1,\dots,a_N\}''$ ($\subset\cM$) is isomorphic
to the von Neumann algebra $\pi_\mu\bigl(\cA_R^{(N)}\bigr)''$ by the
isomorphism defined by
$a_i\mapsto\pi_\mu(X_i)$, $1\le i\le N$, where $\pi_\mu$ is the GNS
representation of $\cA_R^{(N)}$ associated with $\mu$. (A related result in
connection with Connes' embedding problem is found in \cite[\S9]{Br}.)

The next proposition in the single variable case is a direct consequence of
(II) in Section 1.

\begin{prop}\label{P-4.2}
For every $a\in\cM^{sa}$ with $\|a\|\le R$,
$$
\eta(a)=\eta_R(a)=\chi(a).
$$
\end{prop}

The following is obvious by definition and Proposition \ref{P-3.2}.

\begin{prop}\label{P-4.3}
For each $R>0$, $\eta_R(a_1,\dots,a_N)$ is upper semicontinuous on
$(\cM^{sa})^N$ in strong topology, that is, if $(a_1,\dots,a_N)$ and
$(a_1^{(k)},\dots,a_N^{(k)})$ are in
$(\cM^{sa})^N$ for $k\in\bN$ and $a_i^{(k)}\to a_i$ strongly as $k\to\infty$
for $1\le i\le N$, then
$$
\eta_R(a_1,\dots,a_N)\ge
\limsup_{k\to\infty}\eta_R(a_1^{(k)},\dots,a_N^{(k)}).
$$
\end{prop}

The following properties immediately follow from (4) and (6) of Proposition
\ref{P-2.3}.

\begin{prop}\label{P-4.4}
Let $a_1,\dots,a_N\in\cM^{sa}$.
\begin{itemize}
\item[(1)] If $R>R_1>0$, then
$\eta_{R_1}(a_1,\dots,a_N)\le\eta_R(a_1,\dots,a_N)$.
\item[(2)] For $1\le L<N$ and $R>0$,
\begin{eqnarray*}
\eta_R(a_1,\dots,a_N)
&\le&\eta_R(a_1,\dots,a_L)+\eta_R(a_{L+1},\dots,a_N), \\
\eta(a_1,\dots,a_N)
&\le&\eta(a_1,\dots,a_L)+\eta(a_{L+1},\dots,a_N).
\end{eqnarray*}
\end{itemize}
\end{prop}

For Voiculescu's free entropy, it is known that
$\chi_R(a_1,\dots,a_N)=\chi(a_1,\dots,a_N)$ whenever $R\ge\max_i\|a_i\|$.
However, the similar property for $\eta_R$ that
$\eta_R(a_1,\dots,a_N)=\eta(a_1,\dots,a_N)$ for such $R$ is unknown.

\begin{thm}\label{T-4.5} \
\begin{itemize}
\item[(1)] For every $a_1,\dots,a_N\in\cM^{sa}$ and $R>0$,
\begin{eqnarray*}
\eta_R(a_1,\dots,a_N)&\ge&\chi_R(a_1,\dots,a_N), \\
\eta(a_1,\dots,a_N)&\ge&\chi(a_1,\dots,a_N).
\end{eqnarray*}
\item[(2)] If $a_1,\dots,a_N\in\cM^{sa}$ are free and $R\ge\max_i\|a_i\|$,
then
$$
\eta(a_1,\dots,a_N)=\eta_R(a_1,\dots,a_N)=\chi(a_1,\dots,a_N).
$$
\end{itemize}
\end{thm}

\proof
(1)\enspace Let $p\in\bC\<X_1,\dots,X_N\>^{sa}$ and $\delta>0$ be given. One
can choose $r\in\bN$ and $\eps>0$ such that, for each $n\in\bN$,
$(A_1,\dots,A_N)\in\Gamma_R(a_1,\dots,a_N;n,r,\eps)$ implies
$$
|\tr_n(p(A_1,\dots,A_N))-\tau(p(a_1,\dots,a_N))|<\delta.
$$
This gives
\begin{eqnarray*}
&&\int_{(M_n^{sa})_R^N}\exp\bigl(-n^2
\tr_n(p(A_1,\dots,A_N))\bigr)\,d\Lambda_n^{\otimes N} \\
&&\qquad\ge\exp\bigl(-n^2\tau(p(a_1,\dots,a_N))-n^2\delta\bigr)
\,\Lambda_n^{\otimes N}\bigl(\Gamma_R(a_1,\dots,a_N;n,r,\eps)\bigr)
\end{eqnarray*}
so that
\begin{eqnarray*}
&&{1\over n^2}\log\int_{(M_n^{sa})_R^N}\exp\bigl(-n^2
\tr_n(p(A_1,\dots,A_N))\bigr)\,d\Lambda_n^{\otimes N} \\
&&\qquad\ge-\mu_{(a_1,\dots,a_N)}(p)-\delta+{1\over n^2}\log
\Lambda_n^{\otimes N}\bigl(\Gamma_R(a_1,\dots,a_N;n,r,\eps)\bigr).
\end{eqnarray*}
Therefore,
$$
\mu_{(a_1,\dots,a_N)}(p)+\pi_R(p)+\delta\ge\chi_R(a_1,\dots,a_N).
$$
Since $p\in\bC\<X_1,\dots,X_N\>^{sa}$ and $\delta>0$ are arbitrary, we have
$$
\eta_R(a_1,\dots,a_N)\ge\chi_R(a_1,\dots,a_N),
$$
which gives the other inequality as well.

(2)\enspace We have
\begin{eqnarray*}
\eta_R(a_1,\dots,a_N)
&\le&\eta_R(a_1)+\dots+\eta_R(a_N)
\quad\mbox{(by Proposition \ref{P-4.4}\,(2))} \\
&=&\chi(a_1)+\dots+\chi(a_N)
\qquad\mbox{(by Proposition \ref{P-4.2})} \\
&=&\chi(a_1,\dots,a_N)
\end{eqnarray*}
by the additivity of $\chi(a_1,\dots,a_N)$ in the free case (\cite{V2}).
The converse inequality is in (1).\qed

\bigskip
By Proposition \ref{P-4.4} and Theorem \ref{T-4.5} we notice that
$\eta(a_1,\dots,a_N)$ admits the same maximal value as $\chi(a_1,\dots,a_N)$
when restricted on $\|a_i\|\le R$. In fact, the maximum is attained when
$a_1,\dots,a_N$ are free and each $a_i$ has the arcsine distribution
supported on $[-R,R]$.

\begin{remark}\label{R-4.6}{\rm
For $\mu\in\cT\bigl(\cA_R^{(N)}\bigr)$ we denote by $\chi(\mu)$ the free
entropy $\chi(\pi_\mu(X_1),\dots,\pi_\mu(X_N))$ via the GNS representation
$\pi_\mu$ of $\cA_R^{(N)}$ associated with $\mu$. According to \cite{V3},
if $\chi(a_1,\dots,a_N)>-\infty$ for selfadjoint variables $a_1,\dots,a_N$
in a tracial $W^*$-probability space, then $\{a_1,\dots,a_N\}''$ is a
(non-hyperfinite, even non-$\Gamma$) II$_1$ factor. This shows in particular
that if $\mu\in\cT\bigl(\cA_R^{(N)}\bigr)$ is not factorial (or not extremal
in $\cT\bigl(\cA_R^{(N)}\bigr)$), then $\chi(\mu)=-\infty$. Choose two
different $\mu_1,\mu_2\in\cT\bigl(\cA_R^{(N)}\bigr)$ such that $\chi(\mu_1)$
and $\chi(\mu_2)$ are finite; then for $\mu_0:=(\mu_1+\mu_2)/2$ we get
$\eta(\mu_0)>-\infty$ thanks to the concavity of $\eta$ (Theorem
\ref{T-3.4}\,(1)), but $\chi(\mu_0)=-\infty$. Hence, $\eta$ and $\chi$ are
not equal in general; $\chi(\mu)$ for $\mu\in\cT\bigl(\cA_R^{(N)}\bigr)$ is
far from concave.
}\end{remark}

The following two propositions are analogous to \cite[Propositions 3.6 and
3.8]{V2}. But it does not seem that $\eta$ enjoys the change of variable
formulas established in \cite{V2,V4} for $\chi$.

\begin{prop}\label{P-4.7} \
\begin{itemize}
\item[(1)] If $A=[\alpha_{ij}]_{i,j=1}^N\in M_N(\bR)$ and
$\beta_1,\dots,\beta_N\in\bR$, then
$$
\eta\Biggl(\sum_{j=1}^N\alpha_{1j}a_j+\beta_1\1,\dots,
\sum_{j=1}^N\alpha_{Nj}a_j+\beta_N\1\Biggr)
=\eta(a_1,\dots,a_N)+\log|\det A|.
$$
\item[(2)] If $a_1,\dots,a_N$ are linearly independent, then
$\eta(a_1,\dots,a_N)=-\infty$.
\end{itemize}
\end{prop}

\proof
(1)\enspace Write $b_i:=\sum_{j=1}^N\alpha_{ij}a_j+\beta_j\1$ for
$1\le i\le N$. Here, assume that $A$ is invertible (the singular case will
be seen after proving (2)), and define $\Psi:(M_n^{sa})^N\to(M_n^{sa})^N$ by
$$
\Psi(A_1,\dots,A_N)
:=\Biggl(\sum_{j=1}^N\alpha_{1j}A_j+\beta_1I,\dots,
\sum_{j=1}^N\alpha_{Nj}A_j+\beta_NI\Biggr).
$$
It is known (see \cite[6.2.2]{HP1}) that
$d(\Lambda_n^{\otimes N}\circ\Psi)/d\Lambda_n^{\otimes N}=|\det A|^{n^2}$.
For each $R>0$ and $q\in\bC\<X_1,\dots,X_N\>^{sa}$ set
$$
R':=\max_{1\le i\le N}\Biggl(\sum_{j=1}^N|\alpha_{ij}|R+|\beta_i|\Biggr),
$$
$$
p(X_1,\dots,X_N):=q\Biggl(\sum_{j=1}^N\alpha_{1j}X_j+\beta_1\1,
\dots,\sum_{j=1}^N\alpha_{Nj}X_j+\beta_N\1\Biggr).
$$
Then we get
\begin{eqnarray*}
\pi_{R'}(q)&\ge&\limsup_{n\to\infty}\Biggl({1\over n^2}\log
\int_{\Psi\bigl((M_n^{sa})_R^N\bigr)}\exp\bigl(-n^2
\tr_n(q(B_1,\dots,B_N))\bigr)\,d\Lambda_n^{\otimes N}
+{N\over2}\log n\Biggr) \\
&=&\limsup_{n\to\infty}\Biggl({1\over n^2}\log\int_{(M_n^{sa})_R^N}
\exp\bigl(-n^2\tr_n((q\circ\Psi)(A_1,\dots,A_N))
\bigr)\,d(\Lambda_n^{\otimes N}\circ\Psi) \\
&&\hskip8cm+{N\over2}\log n\Biggr) \\
&=&\pi_R(p)+\log|\det A|
\end{eqnarray*}
so that
\begin{eqnarray*}
\tau(q(b_1,\dots,b_N))+\pi_{R'}(q)
&\ge&\tau(p(a_1,\dots,a_N))+\pi_R(p)+\log|\det A| \\
&\ge&\eta_R(a_1,\dots,a_N)+\log|\det A|.
\end{eqnarray*}
Therefore,
$$
\eta_{R'}(b_1,\dots,b_N)\ge\eta_R(a_1,\dots,a_N)+\log|\det A|.
$$
By reversing the roles of $(a_1,\dots,a_N)$ and $(b_1,\dots,b_N)$ we have
$$
\eta(b_1,\dots,b_N)=\eta(a_1,\dots,a_N)+\log|\det A|.
$$

(2)\enspace We may assume that $a_1=\alpha_2a_2+\dots+\alpha_Na_N$ with
$\alpha_i\in\bR$. For every $\eps>0$, since
$a_1=\eps a_1+(1-\eps)\alpha_2a_2+\dots+(1-\eps)\alpha_Na_N$, we get by (1)
$$
\eta(a_1,\dots,a_N)=\eta(a_1,\dots,a_N)+\log\eps,
$$
implying $\eta(a_1,\dots,a_N)=-\infty$. Finally, when $A$ is singular, it
follows from (2) just proven that both sides of the equality in (1) are
$-\infty$.\qed

\begin{prop}\label{P-4.8}
Let $a_1,\dots,a_N,b_1,\dots,b_N\in\cM^{sa}$ be such that $b_1=a_1$ and
$b_i-a_i\in\{a_1,\dots,a_{i-1}\}''$ for $2\le i\le N$, and let
$R_0:=\max_i\|b_i-a_i\|$. Then
\begin{eqnarray*}
\eta_R(a_1,\dots,a_N)&\le&\eta_{R+R_0}(b_1,\dots,b_N)
\quad\mbox{if $R\ge\max_i\|a_i\|$}, \\
\eta_R(b_1,\dots,b_N)&\le&\eta_{R+R_0}(a_1,\dots,a_N)
\quad\mbox{if $R\ge\max_i\|b_i\|$}.
\end{eqnarray*}
Hence,
$$
\eta(a_1,\dots,a_N)=\eta(b_1,\dots,b_N).
$$
\end{prop}

\proof
When $R\ge\max_i\|a_i\|$, as noted just after Definition \ref{D-4.1},
$\mu:=\mu_{(a_1,\dots,a_N)}\in\cT\bigl(\cA_R^{(N)}\bigr)$ and
$\{a_1,\dots,a_N\}''$ is isomorphic to $\pi_\mu\bigl(\cA_R^{(N)}\bigr)''$
by the isomorphism $\psi$ defined by $\psi(a_i)=\pi_\mu(X_i)$,
$1\le i\le N$. For each $2\le i\le N$, since
$\psi(b_i-a_i)\in\{\pi_\mu(X_1),\dots,\pi_\mu(X_{i-1})\}''$, using
Kaplansky density theorem and an argument with functional calculus, one can
choose a sequence of selfadjoint elements $h_i^{(k)}$, $k\in\bN$, in
$C^*(X_1,\dots,X_{i-1})$ ($\subset\cA_R^{(N)}$) such that
$\|h_i^{(k)}\|_R\le\|\psi(b_i-a_i)\|=\|b_i-a_i\|$ and
$\pi_\mu(h_i^{(k)})\to\psi(b_i-a_i)$ strongly as $k\to\infty$. Hence,
there exists a sequence $p_i^{(k)}\in\bC\<X_1,\dots,X_{i-1}\>^{sa}$,
$k\in\bN$, such that $\|p_i^{(k)}\|_R\le\|b_i-a_i\|$ and
$\psi\bigl(p_i^{(k)}(a_1,\dots,a_{i-1})\bigr)
=\pi_\mu\bigl(p_i^{(k)}(X_1,\dots,X_{i-1})\bigr)\to\psi(b_i-a_i)$
strongly; hence $p_i^{(k)}(a_1,\dots,a_{i-1})\to b_i-a_i$ strongly as
$k\to\infty$. Set $b_1^{(k)}:=b_1=a_1$ and
$b_i^{(k)}:=a_i+p_i^{(k)}(a_1,\dots,a_{i-1})$ for $2\le i\le N$. Then
$\|b_i^{(k)}\|\le\|a_i\|+\|b_i-a_i\|\le R+R_0$ and $b_i^{(k)}\to b_i$
strongly as $k\to\infty$. Define $\Psi^{(k)}:(M_n^{sa})^N\to(M_n^{sa})^N$,
$\Psi^{(k)}(A_1,\dots,A_N)=(B_1,\dots,B_N)$, by $B_1:=A_1$ and
$B_i:=A_i+p_i^{(k)}(A_1,\dots,A_{i-1})$ for $2\le i\le N$. Notice that
$\Lambda_n^{\otimes N}$ is $\Psi^{(k)}$-invariant:
$\Lambda_n^{\otimes N}\circ\Psi^{(k)}=\Lambda_n^{\otimes N}$. If
$(A_1,\dots,A_N)\in(M_n^{sa})_R^N$, then
$\Psi^{(k)}(A_1,\dots,A_N)\in(M_n^{sa})_{R+R_0}^N$ because of
$\|B_i\|\le\|A_i\|+\|p_i^{(k)}\|_R\le R+R_0$.

For each $q\in\bC\<X_1,\dots,X_N\>^{sa}$, setting
$p^{(k)}\in\bC\<X_1,\dots,X_N\>^{sa}$ by
$$
p^{(k)}(X_1,\dots,X_N):=q\bigl(X_1,X_2+p_2^{(k)}(X_1),
\dots,X_N+p_N^{(k)}(X_1,\dots,X_{N-1})\bigr),
$$
we get
\begin{eqnarray*}
&&\pi_{R+R_0}(q) \\
&&\quad\ge\limsup_{n\to\infty}
\Biggl({1\over n^2}\log\int_{\Psi\bigl((M_n^{sa})_R^N\bigr)}\exp\bigl(-n^2
\tr_n(q(B_1,\dots,B_N))\bigr)\,d\Lambda_n^{\otimes N}
+{N\over2}\log n\Biggr) \\
&&\quad=\limsup_{n\to\infty}\Biggl({1\over n^2}\int_{(M_n^{sa})_R^N}
\exp\bigl(-n^2\tr_n((q\circ\Psi^{(k)})(A_1,\dots,A_N))\bigr)
\,d\Lambda_n^{\otimes N}+{N\over2}\log n\Biggr) \\
&&\quad=\limsup_{n\to\infty}\Biggl({1\over n^2}\int_{(M_n^{sa})_R^N}
\exp\bigl(-n^2\tr_n(p^{(k)}(A_1,\dots,A_N))\bigr)
\,d\Lambda_n^{\otimes N}+{N\over2}\log n\Biggr) \\
&&\quad=\pi_R(p^{(k)}).
\end{eqnarray*}
Therefore, since $p^{(k)}(a_1,\dots,a_N)=q(b_1^{(k)},\dots,b_N^{(k)})$,
\begin{eqnarray*}
\tau\bigl(q(b_1^{(k)},\dots,b_N^{(k)})\bigr)+\pi_{R+R_0}(q)
&\ge&\tau\bigl(p^{(k)}(a_1,\dots,a_N)\bigr)+\pi_R(p^{(k)}) \\
&\ge&\eta_R(a_1,\dots,a_N)
\end{eqnarray*}
so that letting $k\to\infty$ gives
$$
\tau(q(b_1,\dots,b_N))+\pi_{R+R_0}(q)\ge\eta_R(a_1,\dots,a_N),
$$
which yields
$$
\eta_{R+R_0}(b_1,\dots,b_N)\ge\eta_R(a_1,\dots,a_N).
$$

Since the assumption implies that $a_i-b_i\in\{b_1,\dots,b_{i-1}\}''$ for
$2\le i\le N$, the other direction follows by reversing the roles of
$(a_1,\dots,a_N)$ and $(b_1,\dots,b_N)$.\qed

\section{Equality $\eta=\chi$ for $R$-diagonal variables}
\setcounter{equation}{0}

Let $(\cM,\tau)$ be a tracial $W^*$-probability space as before.
When $x$ is a non-selfadjoint element in $\cM$, Voiculescu's (microstate)
free entropy of $x$ is defined as $\chi(a_1,a_2)$ where $a_1$ and $a_2$
are the real and imaginary parts of $x$, i.e., $a_1=(x+x^*)/2$ and
$a_2=(x-x^*)/2i$. When $a_1$ and $a_2$ are free, we have
$\eta(a_1,a_2)=\chi(a_1,a_2)$ by Theorem \ref{T-4.5}\,(2). In this section
we show that this equality remains true for {\it $R$-diagonal} variables.
The notion of $R$-diagonal was introduced by A.~Nica and R.~Speicher
\cite{NS}, where it was shown that for $x\in\cM$ with ${\rm ker}\,x=\{0\}$,
$x$ is $R$-diagonal if and only if it admits the polar decomposition
$x=u|x|$ with a Haar unitary $u$ free from $|x|=(x^*x)^{1/2}$.

The ``Lebesgue" measure on $M_n$ is defined as
$$
d\hat\Lambda_n(X):=\prod_{i,j=1}^nd(\Re X_{ij})\,d(\Im X_{ij}).
$$
Recall the following correspondences of measures under the transformations
of Descartes decomposition and of polar decomposition (see \cite[6.5.4 and
4.4.7]{HP1}).
\begin{itemize}
\item[(a)] Under the transformation $X\in
M_n\mapsto(A_1,A_2)\in(M_n^{sa})^2$ with $X=A_1+iA_2$, the measure
$\hat\Lambda_n$ on $M_n$ corresponds to the product measure
$$
\Lambda_n\otimes\Lambda_n\quad\mbox{on $(M_n^{sa})^2$}.
$$
\item[(b)] Let $\cU_n$ be the unitary group of order $n$ and $M_n^+$ the
set of positive semidefinite $n\times n$ matrices. Under the transformation
$X\in M_n\mapsto(U,X^*X)\in\cU_n\times M_n^+$ where $U$ is the unitary part
of $X$ (uniquely determined for non-singular matrices $X$, the other case
being $\hat\Lambda_n$-negligible), $\hat\Lambda$ on $M_n$ corresponds to
the product measure
$$
\gamma_n\otimes\bigl(C_n\Lambda_n|_{M_n^+}\bigr)
\quad\mbox{on $\cU_n\times M_n^+$}
$$
where $\gamma_n$ is the Haar probability measure on $\cU_n$ and
$$
C_n:={\pi^{n(n+1)/2}\over2^{n(n-1)/2}\prod_{j=1}^{n-1}j!}.
$$
\end{itemize}

\begin{thm}\label{T-5.1}
For every $x=a_1+ia_2\in\cM$ with $a_1,a_2\in\cM^{sa}$,
$$
\eta(a_1,a_2)\le\chi(x^*x)+{1\over2}\log{\pi\over2}+{3\over4}.
$$
\end{thm}

\proof
For each real polynomial $q(t)$ we set $p\in\bC\<X_1,X_2\>^{sa}$ by
$$
p(X_1,X_2):=q((X_1+iX_2)^*(X_1+iX_2)).
$$
Obviously, we have
\begin{eqnarray*}
q(x^*x)&=&p(a_1,a_2), \\
q(X^*X)&=&p(A_1,A_2)
\quad\mbox{for $X=A_1+iA_2$, $A_1,A_2\in M_n^{sa}$}.
\end{eqnarray*}
For each $R\ge\|x\|$, since
$$
\bigl\{A_1+iA_2:(A_1,A_2)\in(M_n^{sa})_R\bigr\}
\subset(M_n)_{2R}:=\bigl\{X\in M_n:\|X\|\le2R\bigr\},
$$
we get
\begin{eqnarray*}
&&\int_{(M_n^{sa})_R}\exp\bigl(-n^2
\tr_n(p(A_1,A_2))\bigr)\,d\Lambda_n^{\otimes 2}(A_1,A_2) \\
&&\qquad\le\int_{(M_n)_{2R}}\exp\bigl(-n^2
\tr_n(q(X^*X))\bigr)\,d\hat\Lambda_n(X)\quad\mbox{(by (a))} \\
&&\qquad=C_n\int_{(M_n^+)_{4R^2}}\exp\bigl(-n^2
\tr_n(q(A))\bigr)\,d\Lambda_n(A)\quad\mbox{(by (b))}\\
&&\qquad\le C_n\int_{(M_n^{sa})_{4R^2}}\exp\bigl(-n^2
\tr_n(q(A))\bigr)\,d\Lambda_n(A),
\end{eqnarray*}
where $(M_n^+)_{4R^2}:=\bigl\{A\in M_n^+:\|A\|\le4R^2\bigr\}$. Therefore,
\begin{eqnarray*}
\pi_R(p)&=&\limsup_{n\to\infty}\Biggl({1\over n^2}
\log\int_{(M_n^{sa})_R}\exp\bigl(-n^2
\tr_n(p(A_1,A_2))\bigr)\,d\Lambda_n^{\otimes 2}(A_1,A_2)
+\log n\Biggr) \\
&\le&\limsup_{n\to\infty}\Biggl({1\over n^2}
\log\int_{(M_n^{sa})_{4R^2}}\exp\bigl(-n^2
\tr_n(q(A))\bigr)\,d\Lambda_n(A)+{1\over2}\log n\Biggr) \\
&&\qquad+\lim_{n\to\infty}\biggl({1\over n^2}\log C_n
+{1\over2}\log n\biggr) \\
&=&\pi_{4R^2}(q)+{1\over2}\log{\pi\over2}+{3\over4}.
\end{eqnarray*}
This implies that
\begin{eqnarray*}
\eta_R(a_1,a_2)&\le&\tau(p(a_1,a_2))+\pi_R(p) \\
&\le&\tau(q(x^*x))+\pi_{4R^2}(q)+{1\over2}\log{\pi\over2}+{3\over4}.
\end{eqnarray*}
Taking the infimum of the above over $q$ yields
\begin{eqnarray*}
\eta_R(a_1,a_2)&\le&\eta_{4R^2}(x^*x)
+{1\over2}\log{\pi\over2}+{3\over4} \\
&=&\chi(x^*x)+{1\over2}\log{\pi\over2}+{3\over4}
\end{eqnarray*}
by Proposition \ref{P-4.2} thanks to $4R^2\ge\|x^*x\|$, completing the
proof.\qed

\begin{cor}\label{C-5.2}
If $x\in\cM$ and $\chi(x^*x)=-\infty$ (in particular, this is the case if
${\rm ker}\,x\ne\{0\}$), then
$$
\eta(a_1,a_2)=\chi(a_1,a_2)=-\infty.
$$
\end{cor}

\proof
Since
$$
\chi(a_1,a_2)\le\chi(x^*x)+{1\over2}\log{\pi\over2}+{3\over4}=-\infty
$$
by \cite[6.6.3]{HP1}, Theorem \ref{T-5.1} gives the conclusion.\qed

\begin{cor}\label{C-5.3}
If $x=a_1+ia_2$ with $a_1,a_2\in\cM^{sa}$ is $R$-diagonal and $R\ge\|x\|$,
then
$$
\eta(a_1,a_2)=\eta_R(a_1,a_2)=\chi(a_1,a_2).
$$
\end{cor}

\proof
By Corollary \ref{C-5.2} we may assume ${\rm ker}\,x=\{0\}$; then $x$ admits
the polar decomposition as mentioned in the first paragraph of this section.
By Theorem \ref{T-5.1},
$$
\eta_R(a_1,a_2)\le\eta(a_1,a_2)
\le\chi(x^*x)+{1\over2}\log{\pi\over2}+{3\over4}.
$$
But, it is known (see \cite[6.6.8]{HP1} and \cite{NSS}) that
$$
\chi(a_1,a_2)=\chi_R(a_1,a_2)
=\chi(x^*x)+{1\over2}\log{\pi\over2}+{3\over4}.
$$
Since $\chi_R(a_1,a_2)\le\eta_R(a_1,a_2)$ by Theorem \ref{T-4.5}\,(1), we
have the conclusion.\qed

\section{Modified quantity $\tilde\eta(a_1,\dots,a_N)$}
\setcounter{equation}{0}

Although the free entropy-like quantity $\eta(a_1,\dots,a_N)$ is different
from $\chi(a_1,\dots,a_N)$ as mentioned in Remark \ref{R-4.6}, it is a
natural one from the viewpoint of ``free variational principle." In this
section we introduce a modified version of $\eta(a_1,\dots,a_N)$, which is
shown to coincide with $\chi(a_1,\dots,a_N)$ in general. For each $R>0$ we
consider the minimal $C^*$-tensor product
$\cA_R^{(N)}\otimes_\min\cA_R^{(N)}$, whose norm is denoted by
$\|\cdot\|_{R,\min}$. The algebraic tensor product
$\bC\<X_1,\dots,X_N\>\otimes\bC\<X_1,\dots,X_N\>$ is naturally considered
as a (dense) $*$-subalgebra of $\cA_R^{(N)}\otimes_\min\cA_R^{(N)}$. Let
$a_1,\dots,a_N$ be selfadjoint operators in $B(\cH)$. For
$q(X_1,\dots,X_N)=\sum c_{i_1\dots i_k,j_1\dots j_l}
X_{i_1}\cdots X_{i_k}\otimes X_{j_1}\cdots X_{j_l}$
in $\bC\<X_1,\dots,X_N\>\otimes\bC\<X_1,\dots,X_N\>$, we define
$$
q(a_1,\dots,a_N):=\sum c_{i_1,\cdots i_k,j_1\cdots j_l}
a_{i_1}\cdots a_{i_k}\otimes a_{j_1}\cdots a_{j_l}
\in\cA_R^{(N)}\otimes_\min\cA_R^{(N)}.
$$
When $\|a_i\|\le R$, since the functional calculus $h(a_1,\dots,a_N)$ for
$h\in\cA_R^{(N)}$ is a $*$-homomorphism of $\cA_R^{(N)}$ into $B(\cH)$ (see
Section 2), one can define the tensor product $*$-homomorphism
$\Phi\otimes\Phi:\cA_R^{(N)}\otimes_\min\cA_R^{(N)}\to
B(\cH)\overline\otimes B(\cH)=B(\cH\otimes\cH)$. We denote
$(\Phi\otimes\Phi)(g)$ by $g(a_1,\dots,a_N)$ for
$g\in\cA_R^{(N)}\otimes_\min\cA_R^{(N)}$. This extends $q(a_1,\dots,a_N)$
above for $q\in\bC\<X_1,\dots,X_N\>\otimes\bC\<X_1,\dots,X_N\>$.

\begin{definition}\label{D-6.1}{\rm
For each $R>0$ and $g\in\bigl(\cA_R^{(N)}\otimes_\min\cA_R^{(N)}\bigr)^{sa}$
define
\begin{eqnarray*}
&&\pi_R^{(2)}(g):=\limsup_{n\to\infty}\Biggl({1\over n^2}
\log\int_{(M_n^{sa})_R^N}\exp\bigl(-n^2
(\tr_n\otimes\tr_n)(g(A_1,\dots,A_N))\bigr)\,d\Lambda_n^{\otimes N} \\
&&\hskip10cm+{N\over2}\log n\Biggr).
\end{eqnarray*}
}\end{definition}

Properties of $\pi_R^{(2)}(g)$ is similar to those of $\pi_R(h)$ in
Proposition \ref{P-2.3}; for example,
\begin{itemize}
\item[(1)] $\pi_R^{(2)}(g)$ is convex on
$\bigl(\cA_R^{(N)}\otimes_\min\cA_R^{(N)}\bigr)^{sa}$,
\item[(2)] $|\pi_R^{(2)}(g_1)-\pi_R^{(2)}(g_2)|\le\|g_1-g_2\|_{R,\min}$
for all $g_1,g_2\in\bigl(\cA_R^{(N)}\otimes_\min\cA_R^{(N)}\bigr)^{sa}$.
\end{itemize}

\begin{definition}\label{D-6.2}{\rm
Let $\Sigma^{(2)}\<X_1,\dots,X_N\>$ denote the set of selfadjoint linear
functionals $\nu$ on $\bC\<X_1,\dots,X_N\>\otimes\bC\<X_1,\dots,X_N\>$ such
that $\nu(\1)=1$. For each $R>0$ and $\nu\in\Sigma^{(2)}\<X_1,\dots,X_N\>$
define
$$
\eta_R^{(2)}(\nu):=\inf\bigl\{\nu(q)+\pi_R^{(2)}(q):
q\in(\bC\<X_1,\dots,X_N\>\otimes\bC\<X_1,\dots,X_N\>)^{sa}\bigr\}.
$$
}\end{definition}

One can see as Lemma \ref{L-3.3} that if
$\nu\in\Sigma^{(2)}\<X_1,\dots,X_N\>$ and $\eta_R^{(2)}(\nu)>-\infty$,
then $\nu\in\cT\bigl(\cA_R^{(N)}\otimes_\min\cA_R^{(N)}\bigr)$, that is,
$\nu$ (uniquely) extends to a tracial state on
$\cA_R^{(N)}\otimes_\min\cA_R^{(N)}$. Hence the essential domain of
$\eta^{(2)}$ is included in
$\cT\bigl(\cA_R^{(N)}\otimes_\min\cA_R^{(N)}\bigr)$. Furthermore, as in
Theorem \ref{T-3.4}, $\pi_R^{(2)}(g)$ for
$g\in\bigl(\cA_R^{(N)}\otimes_\min\cA_R^{(N)}\bigr)^{sa}$ and
$\eta_R^{(2)}(\nu)$ for
$\nu\in\cT\bigl(\cA_R^{(N)}\otimes_\min\cA_R^{(N)}\bigr)$ are the Legendre
transforms of each other with respect to the Banach space duality between
$\bigl(\cA_R^{(N)}\otimes_\min\cA_R^{(N)}\bigr)^{sa}$ and the selfadjoint
part of $\bigl(\cA_R^{(N)}\otimes_\min\cA_R^{(N)}\bigr)^*$.

\begin{definition}\label{D-6.3}{\rm
For each $a_1,\dots,a_N\in\cM^{sa}$ define
$$
\tilde\eta_R(a_1,\dots,a_N)
:=\eta_R^{(2)}(\mu_{(a_1,\dots,a_N)}\otimes\mu_{(a_1,\dots,a_N)}),
$$
where $\mu_{(a_1,\dots,a_N)}\in\Sigma\<X_1,\dots,X_N\>$ was given in
Definition \ref{D-4.1} and
$\mu_{(a_1,\dots,a_N)}\otimes\mu_{(a_1,\dots,a_N)}$ is an element of
$\Sigma^{(2)}\<X_1,\dots,X_N\>$ defined by the algebraic tensor product of
$\mu_{(a_1,\dots,a_N)}$ and itself. Furthermore, define
$$
\tilde\eta(a_1,\dots,a_N):=\sup_{R>0}\tilde\eta_R(a_1,\dots,a_N).
$$
}\end{definition}

Note that if $R\ge\max_i\|a_i\|$, then
$\mu_{(a_1,\dots,a_N)}\in\cT\bigl(\cA_R^{(N)}\bigr)$ and so
$\mu_{(a_1,\dots,a_N)}\otimes\mu_{(a_1,\dots,a_N)}\in
\cT\bigl(\cA_R^{(N)}\otimes_\min\cA_R^{(N)}\bigr)$.

\begin{thm}\label{T-6.4}
For every $a_1,\dots,a_N\in\cM^{sa}$ and $R>0$,
$$
\eta_R(a_1,\dots,a_N)\ge\tilde\eta_R(a_1,\dots,a_N)
=\chi_R(a_1,\dots,a_N),
$$
$$
\eta(a_1,\dots,a_N)\ge\tilde\eta(a_1,\dots,a_N)
=\chi(a_1,\dots,a_N).
$$
\end{thm}

\proof
By Theorem 4.5\,(1) it is enough to prove $\tilde\eta_R=\chi_R$. (A direct
proof of $\eta_R\ge\tilde\eta_R$ is also easy from definition.) For each
$q\in(\bC\<X_1,\dots,X_N\>\otimes\bC\<X_1,\dots,X_N\>)^{sa}$ notice that
$$
(\mu_{(a_1,\dots,a_N)}\otimes\mu_{(a_1,\dots,a_N)})(q)
=(\tau\otimes\tau)(q(a_1,\dots,a_N))
$$
is a polynomial (of at most order $2$) of mixed moments
$\tau(a_{i_1}\cdots a_{i_k})$ with $k\le K$ for some $K\in\bN$. This is same
for $(\tr_n\otimes\tr_n)(q(A_1,\dots,A_N))$ with $A_1,\dots,A_N\in M_n^{sa}$.
Thus, for any $\delta>0$, one can choose $r\in\bN$ and $\eps>0$ such that,
for each $n\in\bN$, if $(A_1,\dots,A_N)\in\Gamma_R(a_1,\dots,a_N;n,r,\eps)$
then
$$
\big|(\tr_n\otimes\tr_n)(q(A_1,\dots,A_N))
-(\tau\otimes\tau)(q(a_1,\dots,a_N))\big|<\delta.
$$
Then, the proof of $\tilde\eta_R(a_1,\dots,a_N)\ge\chi_R(a_1,\dots,a_N)$ is
the same as in the proof of Theorem \ref{T-4.5}\,(1).

To prove the converse inequality, let $\alpha>\chi_R(a_1,\dots,a_N)$ and
$\beta>0$. There exist $r\in\bN$ and $\eps>0$ such that
$$
\limsup_{n\to\infty}\biggl({1\over n^2}\log\Lambda_n^{\otimes N}
\bigl(\Gamma_R(a_1,\dots,a_N;n,r,\eps)\bigr)+{N\over2}\log n\biggr)<\alpha.
$$
Set
$$
I:=\{(i_1,\dots,i_k):1\le i_1,\dots,i_k\le N,\ 1\le k\le r\}
$$
and define $q\in(\bC\<X_1,\dots,X_N\>\otimes\bC\<X_1,\dots,X_N\>)^{sa}$ by
\begin{eqnarray*}
&&q(X_1,\dots,X_N) \\
&&\quad:={\beta\over\eps^2}\sum_{(i_1,\dots,i_k)\in I}
\bigl(X_{i_1}\cdots X_{i_k}-\tau(a_{i_1}\cdots a_{i_k})\1\bigr)\otimes
\bigl(X_{i_1}\cdots X_{i_k}-\tau(a_{i_1}\cdots a_{i_k})\1\bigr)^*.
\end{eqnarray*}
We notice $(\tau\otimes\tau)(q(a_1,\dots a_N))=0$ and
$$
(\tr_n\otimes\tr_n)(q(A_1,\dots,A_N))
={\beta\over\eps^2}\sum_{(i_1,\dots,i_k)\in I}
\big|\tr_n(A_{i_1}\cdots A_{i_k})-\tau(a_{i_1}\cdots a_{i_k})\big|^2
$$
for $A_1,\dots,A_N\in M_n^{sa}$. Since
$$
(\tr_n\otimes\tr_n)(q(A_1,\dots,A_N))\ge\beta
\quad\mbox{if $(A_1,\dots,A_N)\not\in
\Gamma_R(a_1,\dots,a_N;n,r,\eps)$},
$$
we get
\begin{eqnarray*}
&&\int_{(M_n^{sa})_R^N}\exp\bigl(-n^2
(\tr_n\otimes\tr_n)(q(A_1,\dots,A_N))\bigr)\,d\Lambda_n^{\otimes N} \\
&&\qquad\le\Lambda_n^{\otimes N}\bigl(\Gamma_R(a_1,\dots,a_N;n,r,\eps)
\bigr)+e^{-n^2\beta}\Lambda_n^{\otimes N}\bigl((M_n^{sa})_R^N\bigr)
\end{eqnarray*}
so that
\begin{eqnarray*}
&&\Biggl(\int_{(M_n^{sa})_R^N}\exp\bigl(-n^2(\tr_n\otimes\tr_n)
(q(A_1,\dots,A_N))\bigr)\,d\Lambda_n^{\otimes N}\Biggr)^{1/n^2} \\
&&\qquad\le\Bigl(\Lambda_n^{\otimes N}
\bigl(\Gamma_R(a_1,\dots,a_N;n,r,\eps)\bigr)\Bigr)^{1/n^2}
+e^{-\beta}\Bigl(\Lambda_n^{\otimes N}
\bigl((M_n^{sa})_R^N\bigr)\Bigr)^{1/n^2}.
\end{eqnarray*}
Since
$$
\limsup_{n\to\infty}n^{N/2}\Bigl(\Lambda_n^{\otimes N}
\bigl(\Gamma_R(a_1,\dots,a_N;n,r,\eps)\bigr)\Bigr)^{1/n^2}
\le e^\alpha
$$
and
\begin{eqnarray*}
\lim_{n\to\infty}n^{N/2}\Bigl(\Lambda_n^{\otimes N}
\bigl((M_n^{sa})_R^N\bigr)\Bigr)^{1/n^2}
&=&\lim_{n\to\infty}\Bigl(n^{1/2}\Bigl(\Lambda_n
\bigl((M_n^{sa})_R\bigr)\Bigr)^{1/n^2}\Bigr)^N \\
&=&\Bigl(R(\pi/2)^{1/2}e^{3/4}\Bigr)^N
\end{eqnarray*}
(see \cite[\S5.6]{HP1}), we have
\begin{eqnarray*}
&&\limsup_{n\to\infty}n^{N/2}\Biggl(\int_{(M_n^{sa})_R^N}
\exp\bigl(-n^2(\tr_n\otimes\tr_n)(q(A_1,\dots,A_N))\bigr)
\,d\Lambda_n^{\otimes N}\Biggr)^{1/n^2} \\
&&\qquad\le e^\alpha+e^{-\beta}
\Bigl(R(\pi/2)^{1/2}e^{3/4}\Bigr)^N
\end{eqnarray*}
so that
$$
\pi_R^{(2)}(q)\le\alpha+\log\Bigl(1+e^{-\alpha-\beta}
\Bigl(R(\pi/2)^{1/2}e^{3/4}\Bigr)^N\Bigr).
$$
Since $(\tau\otimes\tau)(q(a_1,\dots,a_N))=0$, this implies that
$$
\tilde\eta_R(a_1,\dots,a_N)
\le\alpha+\log\Bigl(1+e^{-\alpha-\beta}
\Bigl(R(\pi/2)^{1/2}e^{3/4}\Bigr)^N\Bigr).
$$
Letting $\beta\to+\infty$ and then $\alpha\searrow\chi_R(a_1,\dots,a_N)$ we
obtain $\tilde\eta_R(a_1,\dots,a_N)\le\chi_R(a_1,\dots,a_N)$, completing the
proof.\qed

\bigskip
Some known properties of $\chi(a_1,\dots,a_N)$ may be shown based on Theorem
\ref{T-6.4}. For instance, the change of variable formulas in \cite{V2,V4}
can be proven for $\tilde\eta(a_1,\dots,a_N)$ in a bit more easily than for
$\chi(a_1,\dots,a_N)$.

\section{Gibbs ensemble asymptotics}
\setcounter{equation}{0}

For each $R>0$ and $n\in\bN$ we define the ``micro" pressure function
$P_{R,n}(h)$ for $h\in\bigl(\cA_R^{(N)}\bigr)^{sa}$ by
$$
P_{R,n}(h):=\log\int_{(M_n^{sa})_R^N}\exp\bigl(-n^2
\tr_n(h(A_1,\dots,A_N))\bigr)\,d\Lambda_n^{\otimes N}(A_1,\dots,A_N).
$$
This is nothing but the usual (classical) pressure of the continuous
function $n^2\tr_n(h(A_1,\dots,A_N))$ on $(M_n^{sa})_R^N$. The definition
\eqref{F-2.1} is
\begin{equation}\label{F-7.1}
\pi_R(h)=\limsup_{n\to\infty}\biggl({1\over n^2}P_{R,n}(h)
+{N\over2}\log n\biggr)
\quad\mbox{for $h\in\bigl(\cA_R^{(N)}\bigr)^{sa}$}.
\end{equation}
For each $h_0\in\bigl(\cA_R^{(N)}\bigr)^{sa}$ we define the Gibbs
probability measure $\lambda_{R,n}^{h_0}$ on $(M_n^{sa})_R^N$ corresponding
to the function $n^2\tr_n(h_0(A_1,\dots,A_N))$ by
\begin{eqnarray*}
d\lambda_{R,n}^{h_0}(A_1,\dots,A_N)
&:=&{1\over Z_{R,n}^{h_0}}\exp\bigl(-n^2
\tr_n(h_0(A_1,\dots,A_N))\bigr) \\
&&\qquad\quad\times\chi_{(M_n^{sa})_R^N}(A_1,\dots,A_N)
\,d\Lambda_n^{\otimes N}(A_1,\dots,A_N),
\end{eqnarray*}
where the normalization constant
$Z_{R,n}^{h_0}:=\exp\bigl(P_{R,n}(h_0)\bigr)$. Furthermore, we define
$$
\mu_{R,n}^{h_0}(h):=\int_{(M_n^{sa})_R^N}
\tr_n(h(A_1,\dots,A_N))\,d\lambda_{R,n}^{h_0}(A_1,\dots,A_N)
\quad\mbox{for $h\in\bigl(\cA_R^{(N)}\bigr)^{sa}$}.
$$
It is immediate to see that $\mu_{R,n}^{h_0}\in\cT\bigl(\cA_R^{(N)}\bigr)$.
The next lemma is elementary and well known.

\begin{lemma}\label{L-7.1}
With the above definitions, the function $P_{R,n}$ is convex on
$\bigl(\cA_R^{(N)}\bigr)^{sa}$ and differentiable at any
$h_0\in\bigl(\cA_R^{(N)}\bigr)^{sa}$. The supporting function of $P_{R,n}$
at $h_0$ is
$$
-n^2\mu_{R,n}^{h_0}(h)+S(\lambda_{R,n}^{h_0}),
\qquad h\in\bigl(\cA_R^{(N)}\bigr)^{sa},
$$
where $S(\lambda_{R,n}^{h_0})$ is the Boltzmann-Gibbs entropy:
$$
S(\lambda_{R,n}^{h_0}):=-\int_{(M_n^{sa})_R^N}
{d\lambda_{R,n}^{h_0}\over d\Lambda_n^{\otimes N}}
\log{d\lambda_{R,n}^{h_0}\over d\Lambda_n^{\otimes N}}
\,d\Lambda_n^{\otimes N}.
$$
\end{lemma}

\begin{prop}\label{P-7.2}
Let $h_0\in\cA_R^{sa}$.
\begin{itemize}
\item[(1)] There exists an equilibrium tracial state
$\mu_0\in\cT\bigl(\cA_R^{(N)}\bigr)$ associated with $h_0$ such that
$$
\liminf_{n\to\infty}\biggl(
{1\over n^2}S(\lambda_{R,n}^{h_0})+{N\over2}\log n\biggr)
\le\eta_R(\mu_0)\le\limsup_{n\to\infty}\biggl(
{1\over n^2}S(\lambda_{R,n}^{h_0})+{N\over2}\log n\biggr).
$$
\item[(2)] Assume that the limit
$$
\pi_R(h_0)=\lim_{n\to\infty}\biggl(
{1\over n^2}P_{R,n}(h_0)+{N\over2}\log n\biggr)
$$
exists. Then any limit point of the sequence
$\bigl\{{1\over n^2}S(\lambda_{R,n}^{h_0})+{N\over2}\log n\bigr\}$ is
attained as the value $\eta_R(\mu_0)$ of some equilibrium tracial state
$\mu_0$ associated with $h_0$. In particular, if $\mu_0$ is a unique
equilibrium tracial state associated with $h_0$, then
\begin{equation}\label{F-7.2}
\eta_R(\mu_0)=\lim_{n\to\infty}\biggl(
{1\over n^2}S(\lambda_{R,n}^{h_0})+{N\over2}\log n\biggr).
\end{equation}
\end{itemize}
\end{prop}

\proof
(1)\enspace First, note that $\cT\bigl(\cA_R^{(N)}\bigr)$ is a compact
metrizable space in the weak* topology. Hence, there exists a subsequence
$\{n(k)\}$ of $\{n\}$ such that
$$
\pi_R(h_0)=\lim_{k\to\infty}\biggl({1\over n(k)^2}
P_{R,n(k)}(h_0)+{N\over2}\log n(k)\biggr)
$$
and $\mu_{R,n(k)}^{h_0}$ weakly* converges to some
$\mu_0\in\cT\bigl(\cA_R^{(N)}\bigr)$. By Lemma \ref{L-7.1} we get
\begin{eqnarray*}
\mu_0(h_0)+\pi_R(h_0)
&=&\lim_{k\to\infty}\biggl(\mu_{R,n(k)}^{h_0}(h_0)
+{1\over n(k)^2}P_{R,n(k)}(h_0)+{N\over2}\log n(k)\biggr) \\
&=&\lim_{k\to\infty}\biggl({1\over n(k)^2}
S(\lambda_{R,n(k)}^{h_0})+{N\over2}\log n(k)\biggr)
\end{eqnarray*}
and for every $h\in\bigl(\cA_R^{(N)}\bigr)^{sa}$
\begin{eqnarray*}
\mu_0(h)+\pi_R(h)
&=&\lim_{k\to\infty}\mu_{R,n(k)}^{h_0}(h)
+\limsup_{n\to\infty}\biggl({1\over n^2}P_{R,n}(h)
+{N\over2}\log n\biggr) \\
&\ge&\limsup_{k\to\infty}\biggl(\mu_{R,n(k)}^{h_0}(h)
+{1\over n(k)^2}P_{R,n(k)}(h)+{N\over2}\log n(k)\biggr) \\
&\ge&\lim_{k\to\infty}\biggl({1\over n(k)^2}
S(\lambda_{R,n(k)}^{h_0})+{N\over2}\log n(k)\biggr)
=\mu_0(h_0)+\pi_R(h_0).
\end{eqnarray*}
By Theorem \ref{T-3.4}\,(1) this implies that $\mu_0$ is an equilibrium
tracial state associated with $h_0$ and
$\eta_R(\mu_0)=\mu_0(h_0)+\pi_R(h_0)$. Hence we have the conclusion.

(2)\enspace Let $\{n(k)\}$ be a subsequence of $\{n\}$ for which the limit
$$
\alpha:=\lim_{k\to\infty}\biggl({1\over n(k)^2}S(\lambda_{R,n(k)}^{h_0})
+{N\over2}\log n(k)\biggr)
$$
exists. We may assume that $\mu_{R,n(k)}^{h_0}\to\mu_0$ weakly* for some
$\mu_0\in\cT\bigl(\cA_R^{(N)}\bigr)$. Then, as in the proof of (1), we get
$\mu_0(h_0)+\pi_R(h_0)=\alpha\le\mu_0(h)+\pi_R(h)$ for all
$h\in\bigl(\cA_R^{(N)}\bigr)^{sa}$, and the result follows.\qed

\bigskip
Although the formula \eqref{F-7.2} has been shown under the existence of
limit in the definition of $\pi_R(h_0)$ as well as the strong assumption
of unique equilibrium, the above proposition provides the noncommutative
multivariate version of the property (IV) in Section 1. The existence of
limit of this kind seems one of the major questions in random matrix
theory, and recently a similar limit behavior has been investigated by
A.~Guionnet \cite{Gu} for several particular noncommuting polynomials of
interest.

To get rid of the convergence problem, we are tempted to introduce the
free pressure by using the limit via a ultrafilter, as Voiculescu defined
the free entropy $\chi^\omega$ in \cite{Va}. Let $\omega$ be a fixed free
ultrafilter, i.e., $\omega\in\beta\bN\setminus\bN$. Since the inside of
the $\limsup$ in \eqref{F-7.1} is a bounded sequence, we can define
$$
\pi_R^\omega(h):=\lim_{n\to\omega}\biggl(
{1\over n^2}P_{R,n}(h)+{N\over2}\log n\biggr)
$$
for each $R>0$ and $h\in\cA_R^{sa}$. Then $\pi_R^\omega$ has the same
properties as $\pi_R$, and we define its Legendre transform with respect to
the duality between $\bigl(\cA_R^{(N)}\bigr)^{sa}$ and
$\bigl(\cA_R^{(N)}\bigr)^{*,sa}$ by
$$
\eta_R^\omega(\mu):=\inf\bigl\{\mu(h)+\pi_R^\omega(h):
h\in\bigl(\cA_R^{(N)}\bigr)^{sa}\bigr\}
\quad\mbox{for $\mu\in\bigl(\cA_R^{(N)}\bigr)^{*,sa}$}.
$$
As before we have $\eta_R^\omega(\mu)=-\infty$ unless
$\mu\in\cT\bigl(\cA_R^{(N)}\bigr)$ and
$$
\pi_R^\omega(h)=\sup\bigl\{-\mu(h)+\eta_R^\omega(\mu):
\mu\in\cT\bigl(\cA_R^{(N)}\bigr)\bigr\}
\quad\mbox{for $h\in\bigl(\cA_R^{(N)}\bigr)^{sa}$}.
$$
Similarly, we define for
$g\in\bigl(\cA_R^{(N)}\otimes_\min\cA_R^{(N)}\bigr)^{sa}$
\begin{eqnarray*}
&&\pi_R^{(2)}(g):=\limsup_{n\to\omega}\Biggl({1\over n^2}
\log\int_{(M_n^{sa})_R^N}\exp\bigl(-n^2
(\tr_n\otimes\tr_n)(g(A_1,\dots,A_N))\bigr)\,d\Lambda_n^{\otimes N} \\
&&\hskip10cm+{N\over2}\log n\Biggr).
\end{eqnarray*}
and for $\nu\in\cT\bigl(\cA_R^{(N)}\otimes_\min\cA_R^{(N)}\bigr)$
$$
\eta_R^{(2),\omega}(\nu):=\inf\bigl\{\nu(h)+\pi_R^{(2),\omega}(h):
h\in(\cA_R\otimes\cA_R^{op})^{sa}\bigr\}.
$$
Furthermore, for each $a_1,\dots,a_N\in\cM^{sa}$ define
\begin{eqnarray*}
\eta_R^\omega(a_1,\dots,a_N)&:=&\eta_R^\omega(\mu_{(a_1,\dots,a_M)}), \\
\eta^\omega(a_1,\dots,a_N)&:=&\sup_{R>0}\eta_R^\omega(a_1,\dots,a_N), \\
\tilde\eta_R^\omega(a_1,\dots,a_N)&:=&
\eta_R^{(2),\omega}(\mu_{(a_1,\dots,a_N)}\otimes\mu_{(a_1,\dots,a_N)}), \\
\tilde\eta^\omega(a_1,\dots,a_N)&:=&
\sup_{R>0}\tilde\eta_R^\omega(a_1,\dots,a_N).
\end{eqnarray*}

The next proposition can be shown as Theorems \ref{T-4.5} and \ref{T-6.4}.

\begin{prop}\label{P-7.3}
For every $a_1,\dots,a_N\in\cM^{sa}$ and $R>0$,
$$
\eta_R^\omega(a_1,\dots,a_N)\ge\tilde\eta_R^\omega(a_1,\dots,a_N)
=\chi_R^\omega(a_1,\dots,a_N),
$$
$$
\eta^\omega(a_1,\dots,a_N)\ge\tilde\eta^\omega(a_1,\dots,a_N)
=\chi^\omega(a_1,\dots,a_N).
$$
\end{prop}

\begin{prop}\label{P-7.4}
Let $h_0\in\bigl(\cA_R^{(N)}\bigr)^{sa}$ and define
$$
\mu_{R,\omega}^{h_0}(h):=\lim_{n\to\omega}\mu_{R,n}^{h_0}(h),
\qquad h\in\bigl(\cA_R^{(N)}\bigr)^{sa}.
$$
Then $\mu_{R,\omega}^{h_0}$ is an equilibrium tracial state associated with
$h_0$ in the sense that
$$
\pi_R^\omega(h_0)=-\mu_{R,\omega}^{h_0}(h_0)
+\eta_R^\omega(\mu_{R,\omega}^{h_0}).
$$
Moreover,
$$
\eta_R^\omega(\mu_{R,\omega}^{h_0})=\lim_{n\to\omega}\biggl(
{1\over n^2}S(\lambda_{R,n}^{h_0})+{N\over2}\log n\biggr).
$$
\end{prop}

\proof
It is immediate to see that
$\mu_{R,\omega}^{h_0}\in\cT\bigl(\cA_R^{(N)}\bigr)$. By Lemma \ref{L-7.1}
we get for every $h\in\bigl(\cA_R^{(N)}\bigr)^{sa}$
\begin{eqnarray*}
\mu_{R,\omega}^{h_0}(h)+\pi_R^\omega(h)
&=&\lim_{n\to\omega}\biggl(\mu_{R,n}^{h_0}(h)+{1\over n^2}
P_{R,n}(h)+{N\over2}\log n\biggr) \\
&\ge&\lim_{n\to\omega}\biggl({1\over n^2}S(\lambda_{R,n}^{h_0})
+{N\over2}\log n\biggr) \\
&=&\mu_{R,\omega}^{h_0}(h_0)+\pi_R^\omega(h_0).
\end{eqnarray*}
This implies the conclusion.\qed

\bigskip
Finally, it is worth noting that the Gibbs ensemble asymptotics (or random
matrix approximation) provides a useful tool to obtain free probabilistic
analogs of classical problems. For example, the free transportation cost
inequality established by Ph.~Biane and D.~Voiculescu \cite{BV} in case of
single variables was re-proven in \cite{HPU} by using this tool, and the
multivariate case will be discussed in our forthcoming paper.

\end{document}